\documentclass[12pt]{article}
\usepackage{amsmath}
\usepackage{amssymb}
\usepackage{t1enc}
\usepackage[latin1]{inputenc}
\usepackage[english]{babel}
\usepackage{amsfonts}
\usepackage{latexsym}
\usepackage{bm}
\newtheorem{theorem}{Theorem}[section]

\newtheorem{lemma}[theorem]{Lemma}

\newtheorem{conj}[theorem]{Conjecture}

\begin{document}
\title{Intersecting families, cross-intersecting families, and a proof of a conjecture of Feghali, Johnson and Thomas}

\author{Peter Borg\\[5mm]
Department of Mathematics \\
University of Malta\\
Malta\\
\texttt{peter.borg@um.edu.mt}}
\date{} \maketitle

\begin{abstract} A family $\mathcal{A}$ of sets is said to be
\emph{intersecting} if every two sets in $\mathcal{A}$ intersect. Two families $\mathcal{A}$ and $\mathcal{B}$ are said to be \emph{cross-intersecting} if each set in $\mathcal{A}$ intersects each set in $\mathcal{B}$. For a positive integer $n$, let $[n] = \{1, \dots, n\}$ and $\mathcal{S}_n = \{A \subseteq [n] \colon 1 \in A\}$. In this note, we extend the Erd\H{o}s--Ko--Rado Theorem by showing that if $\mathcal{A}$ and $\mathcal{B}$ are non-empty cross-intersecting families of subsets of $[n]$, $\mathcal{A}$ is intersecting, and $a_0, a_1, \dots, a_n, b_0, b_1, \dots, b_n$ are non-negative real numbers such that $a_i + b_i \geq a_{n-i} + b_{n-i}$ and $a_{n-i} \geq b_i$ for each $i \leq n/2$, then 
\[\sum_{A \in \mathcal{A}} a_{|A|} + \sum_{B \in \mathcal{B}} b_{|B|} \leq \sum_{A \in \mathcal{S}_n} a_{|A|} + \sum_{B \in \mathcal{S}_n} b_{|B|}.\]
For a graph $G$ and an integer $r$, let ${\mathcal{I}_G}^{(r)}$ denote the family of $r$-element independent sets of $G$. Inspired by a problem of Holroyd and Talbot, Feghali, Johnson and Thomas conjectured that if $r < n$ and $G$ is a \emph{depth-two claw} with $n$ leaves, then $G$ has a vertex $v$ such that $\{A \in {\mathcal{I}_G}^{(r)} \colon v \in A\}$ is a largest intersecting subfamily of ${\mathcal{I}_G}^{(r)}$. They proved this for $r \leq \frac{n+1}{2}$. We use the result above to prove the full conjecture.
\end{abstract}

\section{Introduction}

Unless otherwise stated, we shall use small letters such as $x$ to
denote non-negative integers or elements of a set,
capital letters such as $X$ to denote sets, and calligraphic
letters such as $\mathcal{F}$ to denote \emph{families}
(that is, sets whose members are sets themselves). It is to be
assumed that arbitrary sets and families are finite. We
call a set $A$ an \emph{$r$-element set} if its size $|A|$ is $r$, that is, if it contains exactly $r$ elements (also called members).

The set $\{1, 2, \dots\}$ of positive integers is denoted by $\mathbb{N}$. For any integer $n \geq 0$, the set $\{i \in \mathbb{N}
\colon i \leq n\}$ is denoted by $[n]$. Note that $[0]$ is the empty set $\emptyset$. For a set $X$, the \emph{power set of $X$} (that is, $\{A \colon A \subseteq X\}$) is denoted by $2^X$. The family of $r$-element subsets of $X$ is denoted by $X \choose r$. The family of $r$-element sets in a family $\mathcal{F}$ is denoted by $\mathcal{F}^{(r)}$. If $\mathcal{F} \subseteq 2^X$ and $x \in X$, then we denote the family $\{F \in \mathcal{F} \colon x \in F\}$ by $\mathcal{F}(x)$. We call $\mathcal{F}(x)$ a \emph{star of $\mathcal{F}$} if $\mathcal{F}(x) \neq \emptyset$.

We say that a set $A$ \emph{intersects} a set $B$ if $A$ and $B$ have at least one common element (that is, $A \cap B \neq \emptyset$). A family $\mathcal{A}$ is said to be \emph{intersecting} if for every $A, B \in \mathcal{A}$, $A$ and $B$ intersect. The stars of a family $\mathcal{F}$ (with $|\bigcup_{F \in \mathcal{F}} F| \geq 1$) are the simplest intersecting subfamilies of $\mathcal{F}$. We say that $\mathcal{F}$ has the star property if at least one of the largest intersecting subfamilies of $\mathcal{F}$ is a star of $\mathcal{F}$.

One of the most popular endeavours in extremal set theory is that of determining the size of a largest intersecting subfamily of a given family $\mathcal{F}$. This started in \cite{EKR}, which features the following classical result, known as the Erd\H os-Ko-Rado (EKR) Theorem. 

\begin{theorem}[{EKR Theorem \cite{EKR}}] \label{EKRthm} If $r \leq n/2$ and $\mathcal{A}$ is an intersecting subfamily of ${[n] \choose r}$, then $|\mathcal{A}| \leq {n-1 \choose r-1}$.
\end{theorem} 
This means that ${[n] \choose r}$ has the star property. There are various proofs of the EKR Theorem (see \cite{D,K,Kat}), two of which are particularly short and beautiful: Katona's \cite{K}, which introduced the elegant cycle method, and Daykin's \cite{D}, using the fundamental Kruskal--Katona Theorem \cite{Ka,Kr}. The EKR Theorem gave rise to some of the highlights in extremal set theory \cite{AK1,F_t1,Kat,W} and inspired many results that establish how large a system of sets can be under certain intersection conditions; see \cite{Borg7,DF,F2,F,FT,HST,HT}.

If $\mathcal{A}$ and $\mathcal{B}$ are families such that each set in $\mathcal{A}$ intersects each set in $\mathcal{B}$, then $\mathcal{A}$ and $\mathcal{B}$ are said to be \emph{cross-intersecting}. 

For intersecting subfamilies of a given family $\mathcal{F}$,
the natural question to ask is how large they can be. 
A natural variant of this intersection problem is the problem of maximizing the sum or the product of sizes of cross-intersecting subfamilies (not necessarily distinct or non-empty) of $\mathcal{F}$. This has recently attracted much attention. The relation between the original intersection problem, the sum problem and the product problem is studied in \cite{Borg8}. Solutions have been obtained for various families; most of the known results are referenced in \cite{LMS2,LMS3}, which treat the product problem for families of subsets of $[n]$ of size at most $r$. 

Here we consider the sum problem for the case where at least one of two cross-intersecting families $\mathcal{A}$ and $\mathcal{B}$ of subsets of $[n]$ is an intersecting family. We actually consider a more general setting of \emph{weighted} sets, where each set of size $i$ is assigned two non-negative integers $a_i$ and $b_i$, and the objective is to maximize $\sum_{A \in \mathcal{A}} a_{|A|} + \sum_{B \in \mathcal{B}} b_{|B|}$. Let $\mathcal{S}_n$ denote the star $\{A \subseteq [n] \colon 1 \in A\}$ of $2^{[n]}$. In Section~\ref{Proofs}, we prove the following extension of the EKR Theorem.
\begin{theorem}\label{mainresult} If $\mathcal{A}$ and $\mathcal{B}$ are non-empty cross-intersecting families of subsets of $[n]$, $\mathcal{A}$ is intersecting, and $a_0, a_1, \dots, a_n, b_0, b_1, \dots, b_n$ are non-negative real numbers such that $a_i + b_i \geq a_{n-i} + b_{n-i}$ and $a_{n-i} \geq b_i$ for each $i \leq n/2$, then 
\[\sum_{A \in \mathcal{A}} a_{|A|} + \sum_{B \in \mathcal{B}} b_{|B|} \leq \sum_{A \in \mathcal{S}_n} a_{|A|} + \sum_{B \in \mathcal{S}_n} b_{|B|}.\]
\end{theorem}
The EKR Theorem is obtained by taking $r \leq n/2$, $\mathcal{B} = \mathcal{A} \subseteq {[n] \choose r}$, and $b_i = 0 = a_i - 1$ for each $i \in \{0\} \cup [n]$. 

We use Theorem~\ref{mainresult} to prove a conjecture of Feghali, Johnson and Thomas \cite[Conjecture~2.1]{FJT}. Before stating the conjecture, we need some further definitions and notation.

A \emph{graph} $G$ is a pair $(X,\mathcal{Y})$, where $X$ is a set, called the \emph{vertex set of $G$}, and $\mathcal{Y}$ is a subset of ${X \choose 2}$ and is called the \emph{edge set of $G$}. The vertex set of $G$ and the edge set of $G$ are denoted by $V(G)$ and $E(G)$, respectively. An element of $V(G)$ is called a \emph{vertex of $G$}, and an element of $E(G)$ is called an \emph{edge of $G$}. We may represent an edge $\{v,w\}$ by $vw$. If $vw$ is an edge of $G$, then we say that $v$ is \emph{adjacent} to $w$ (in $G$). 
A subset $I$ of $V(G)$ is an \emph{independent set of $G$} if $vw \notin E(G)$ for every $v, w \in I$. Let $\mathcal{I}_G$ denote the family of independent sets of $G$. An independent set $J$ of $G$ is \emph{maximal} if $J \nsubseteq I$ for each independent set $I$ of $G$ such that $I \neq J$. The size of a smallest maximal independent set of $G$ is denoted by $\mu(G)$.

Holroyd and Talbot introduced the problem of determining whether ${\mathcal{I}_G}^{(r)}$ has the star property for a given graph $G$ and an integer $r \geq 1$. The Holroyd--Talbot (HT) Conjecture \cite[Conjecture~7]{HT} claims that ${\mathcal{I}_G}^{(r)}$ has the star property if $\mu(G) \geq 2r$. It is proved in \cite{LMS1} that the conjecture is true if $\mu(G)$ is sufficiently large depending on $r$. 
By the EKR Theorem, it is true if $G$ has no edges. The HT Conjecture has been verified for several classes of graphs \cite{BH1,BH,HHS,HS,HST,HT,HK,T,Woodroofe}. 
As demonstrated in \cite{BH}, for $r > \mu(G)/2$, whether ${\mathcal{I}_G}^{(r)}$ has the star property or not depends on $G$ and $r$ (both cases are possible).

A \emph{depth-two claw} is a graph consisting of $n$ pairwise disjoint edges $x_1y_1, \dots, x_ny_n$ together with a vertex $x_0 \notin \{x_1, \dots, x_n, y_1, \dots, y_n\}$ that is adjacent to each of $y_1, \dots, y_n$. This graph will be denoted by $T_n$. Thus, $T_n = (\{x_0, x_1, \dots, x_n, y_1, \dots, y_n\}, \{x_0y_1, \dots, x_0y_n, x_1y_1, \dots, x_ny_n\})$. For each $i \in [n]$, we may take $x_i$ and $y_i$ to be $(i,1)$ and $(i,2)$, respectively. Let $X_n = \{x_i \colon i \in [n]\}$ and 
\[\mathcal{L}_{n,k} = \left\{\{(i_1,j_1), \dots, (i_r,j_r)\} \colon r \in [n], \{i_1, \dots, i_r\} \in {[n] \choose r}, \, j_1, \dots, j_r \in [k]\right\}.\]
Note that 
\begin{equation}
{\mathcal{I}_{T_n}}^{(r)} = {\mathcal{L}_{n,2}}^{(r)} \cup \left\{A \cup \{x_0\} \colon A \in {X_n \choose r-1}\right\}. \label{ITnr}
\end{equation}

The family ${\mathcal{I}_{T_n}}^{(r)}$ is empty for $r > n+1$, and consists only of the set $\{x_0, x_1, \dots, x_n\}$ for $r = n+1$. In \cite{FJT}, Feghali, Johnson and Thomas showed that ${\mathcal{I}_{T_n}}^{(r)}$ does not have the star property for $r = n$, and they made the following conjecture.

\begin{conj}[{\cite{FJT}}] \label{FJTconj} If $r \leq n-1$, then ${\mathcal{I}_{T_n}}^{(r)}$ has the star property.
\end{conj}
They proved the conjecture for $r \leq \frac{n+1}{2}$.

\begin{theorem}[{\cite{FJT}}] \label{FJTresult} If $n \geq 2r-1$, then
${\mathcal{I}_{T_n}}^{(r)}$ has the star property. 
\end{theorem}
In the next section, we settle the full conjecture using Theorem~\ref{mainresult} for $r > \frac{n+1}{2}$.
\begin{theorem}\label{conjtrue} Conjecture~\ref{FJTconj} is true.
\end{theorem}
Our proof for $r \leq \frac{n}{2}+1$ provides an alternative proof of Theorem~\ref{FJTresult}.

\section{Proofs} \label{Proofs}

In this section, we prove Theorems~\ref{mainresult} and \ref{conjtrue}.\\
\\
\textbf{Proof of Theorem~\ref{mainresult}.} For each $i \in \{0\} \cup [n]$, let $c_i = |\mathcal{A}^{(i)}|a_i + |\mathcal{B}^{(i)}|b_i$. Since $\mathcal{A}$ and $\mathcal{B}$ are non-empty and cross-intersecting, we have $\emptyset \notin \mathcal{A}$ and $\emptyset \notin \mathcal{B}$, so $\mathcal{A}^{(0)} = \mathcal{B}^{(0)} = \emptyset$. Thus, $c_0 = 0$. The result is immediate if $n = 1$. Suppose $n > 1$. 

Consider any positive integer $r \leq n/2$. Since $\mathcal{A}$ and $\mathcal{B}$ are cross-intersecting and $\mathcal{A}$ is intersecting, $[n] \backslash A \notin \mathcal{A}^{(n-r)}$ for each $A \in \mathcal{A}^{(r)} \cup \mathcal{B}^{(r)}$, so 
\begin{align} |\mathcal{A}^{(n-r)}| &\leq {n \choose n-r} - |\{[n] \backslash A \colon A \in \mathcal{A}^{(r)} \cup \mathcal{B}^{(r)}\}| = {n \choose r} - |\mathcal{A}^{(r)} \cup \mathcal{B}^{(r)}| \nonumber \\
&= {n \choose r} - |\mathcal{A}^{(r)}| - |\mathcal{B}^{(r)} \backslash \mathcal{A}^{(r)}|. \nonumber
\end{align}
Similarly, since $\mathcal{A}$ and $\mathcal{B}$ are cross-intersecting, $|\mathcal{B}^{(n-r)}| \leq {n \choose r} - |\mathcal{A}^{(r)}|$. Thus, 
\begin{align} c_r + c_{n-r} &\leq |\mathcal{A}^{(r)}|a_r + |\mathcal{A}^{(r)}|b_r + |\mathcal{B}^{(r)} \backslash \mathcal{A}^{(r)}|b_r \nonumber \\
&\quad + \left({n \choose r} - |\mathcal{A}^{(r)}| - |\mathcal{B}^{(r)} \backslash \mathcal{A}^{(r)}|\right)a_{n-r} + \left({n \choose r} - |\mathcal{A}^{(r)}|\right)b_{n-r} \nonumber \\
&= |\mathcal{A}^{(r)}|(a_r + b_r - a_{n-r} - b_{n-r}) - |\mathcal{B}^{(r)} \backslash \mathcal{A}^{(r)}| (a_{n-r} - b_r) + {n \choose r}(a_{n-r} + b_{n-r}) \nonumber \\
&\leq {n-1 \choose r-1}(a_r + b_r - a_{n-r} - b_{n-r}) + {n \choose r}(a_{n-r} + b_{n-r}) \nonumber
\end{align}
by Theorem~\ref{EKRthm} and the given conditions $a_r + b_r \geq a_{n-r} + b_{n-r}$ and $a_{n-r} \geq b_r$. Therefore, $c_r + c_{n-r} \leq {n-1 \choose r-1}(a_r + b_r) + {n-1 \choose r}(a_{n-r} + b_{n-r}) = {n-1 \choose r-1}(a_r + b_r) + {n-1 \choose n-r-1}(a_{n-r} + b_{n-r})$. Note that if $r = n/2$, then we have $c_{n/2} + c_{n/2} \leq {n-1 \choose n/2-1}(a_{n/2} + b_{n/2}) + {n-1 \choose n/2-1}(a_{n/2} + b_{n/2})$, and hence $c_{n/2} \leq {n-1 \choose n/2-1}(a_{n/2} + b_{n/2})$.

If $n$ is odd, then
\begin{align} \sum_{A \in \mathcal{A}} a_{|A|} + \sum_{B \in \mathcal{B}} b_{|B|} &= \sum_{i=0}^\frac{n-1}{2} (|\mathcal{A}^{(i)}|a_i + |\mathcal{A}^{(n-i)}|a_{n-i} + |\mathcal{B}^{(i)}|b_i + |\mathcal{B}^{(n-i)}|b_{n-i}) \nonumber \\
&= c_0 + c_n + \sum_{i=1}^\frac{n-1}{2} (c_i + c_{n-i}) \nonumber \\
&\leq a_n + b_n + \sum_{i=1}^\frac{n-1}{2} \left({n-1 \choose i-1}(a_i + b_i) + {n-1 \choose n-i-1}(a_{n-i} + b_{n-i})\right) \nonumber \\
&= \sum_{A \in \mathcal{S}_n} a_{|A|} + \sum_{B \in \mathcal{S}_n} b_{|B|}. \nonumber
\end{align}
Similarly, if $n$ is even, then
\begin{align} \sum_{A \in \mathcal{A}} a_{|A|} + \sum_{B \in \mathcal{B}} b_{|B|} &= c_0 + c_n + c_{n/2} + \sum_{i=1}^{\frac{n}{2}-1} (c_i + c_{n-i}) \nonumber \\
&\leq a_n + b_n + {n-1 \choose n/2-1}(a_{n/2} + b_{n/2}) \nonumber \\
& \quad + \sum_{i=1}^{\frac{n}{2}-1} \left({n-1 \choose i-1}(a_i + b_i) + {n-1 \choose n-i-1}(a_{n-i} + b_{n-i})\right) \nonumber \\
&= \sum_{A \in \mathcal{S}_n} a_{|A|} + \sum_{B \in \mathcal{S}_n} b_{|B|}, \nonumber
\end{align}
as required.~\hfill{$\Box$}\\

For $(i,j) \in [n] \times [2,k]$, let $\delta_{i,j} \colon \mathcal{L}_{n,k} \rightarrow \mathcal{L}_{n,k}$ be defined by
\[ \delta_{i,j}(A) = \left\{ \begin{array}{ll}
(A \backslash \{(i,j)\}) \cup \{(i,1)\} & \mbox{if $(i,j) \in A$};\\
A & \mbox{otherwise,}
\end{array} \right. \]
and let $\Delta_{i,j} \colon
2^{\mathcal{L}_{n,k}} \rightarrow 2^{\mathcal{L}_{n,k}}$ be the compression operation defined by
$$\Delta_{i,j}(\mathcal{A}) = \{\delta_{i,j}(A) \colon A \in \mathcal{A}\}
\cup \{A \in \mathcal{A} \colon \delta_{i,j}(A) \in
\mathcal{A}\}.$$
It is well known that $|\Delta_{i,j}(\mathcal{A})| = |\mathcal{A}|$ and  that, if $\mathcal{A}$ is intersecting, then $\Delta_{i,j}(\mathcal{A})$ is intersecting. Moreover, we have the following special case of \cite[Corollary~3.2]{Borg2}.

\begin{lemma}\label{complemma} If $\mathcal{A}$ is an intersecting subfamily of $\mathcal{L}_{n,k}$ and
\[\mathcal{A}^* = \Delta_{n,k} \circ \dots \circ \Delta_{n,2} \circ
\dots \circ \Delta_{1,k} \circ \dots \circ
\Delta_{1,2}(\mathcal{A}),\]
then $|A \cap B \cap X_n| \geq 1$ for any $A, B \in
\mathcal{A}^*$.
\end{lemma}
\textbf{Proof of Theorem~\ref{conjtrue}.} The result is trivial for $0 \leq r \leq 1$, so consider $2 \leq r \leq n-1$. Let $\mathcal{R} = {\mathcal{I}_{T_n}}^{(r)}$. Let $\mathcal{E}$ be an intersecting subfamily of $\mathcal{R}$. Let $\mathcal{F} = \{A \in \mathcal{R} \colon x_1 \in A\}$. For any $\mathcal{H} \subseteq \mathcal{R}$, let $\mathcal{H}_0 = \{H \in \mathcal{H} \colon x_0 \notin H\}$, $\mathcal{H}_1 = \{H \in \mathcal{H} \colon x_0 \in H\}$ and $\mathcal{H}_1' = \{H \backslash \{x_0\} \colon H \in \mathcal{H}_1\}$; by (\ref{ITnr}), $\mathcal{H}_0 \subseteq {\mathcal{L}_{n,2}}^{(r)}$ and $\mathcal{H}_1' \subseteq {X_n \choose r-1}$.\medskip

\textit{Case 1: $n \geq 2r-2$.} Then $r \leq \frac{n}{2} + 1$. Let $\gamma \colon \mathcal{R} \rightarrow \mathcal{R}$ be defined by \[\gamma(A) := \left\{\begin{array}{ll}
(A \backslash \{x_0\}) \cup \{y_1\} & \mbox{if $x_0 \in A, y_1 \notin A$ and $(A \backslash \{x_0\}) \cup \{y_1\} \in \mathcal{R}$;}\\
A & \mbox{otherwise,} \end{array} \right. \]
and let $\Gamma \colon 2^{\mathcal{R}} \rightarrow 2^{\mathcal{R}}$ be the compression operation defined by
\[\Gamma(\mathcal{A}) := \{\gamma(A) \colon A \in
\mathcal{A}\} \cup \{A\in \mathcal{A} \colon \gamma(A) \in
\mathcal{A}\}.\]
Note that $\gamma(A) \neq A$ if and only if $x_0 \in A$ and $x_1, y_1 \notin A$. Let $\mathcal{G} = \Gamma(\mathcal{E})$. Then $\mathcal{G} \subseteq \mathcal{R}$ and $|\mathcal{G}|  = |\mathcal{E}|$. By \cite[Lemma~2.1]{BH}, $\mathcal{G}_0$ and $\mathcal{G}_1'$ are intersecting. Since $r-1 \leq \frac{n}{2}$ and $\mathcal{G}_1' \subseteq {X_n \choose r-1}$, Theorem~\ref{EKRthm} gives us $|\mathcal{G}_1'| \leq |\{A \in {X_n \choose r-1} \colon x_1 \in A\}| = |\mathcal{F}_1'|$. It is well known that, for any $k \geq 2$, ${\mathcal{L}_{n,k}}^{(r)}$ has the star property \cite[Theorem~5.2]{DF}. 
Thus, $|\mathcal{G}_0| \leq |\{A \in {\mathcal{L}_{n,2}}^{(r)} \colon x_1 \in A\}| = |\mathcal{F}_0|$. Therefore, we have $|\mathcal{E}| = |\mathcal{G}| = |\mathcal{G}_0| + |\mathcal{G}_1'| \leq |\mathcal{F}_0| + |\mathcal{F}_1'| = |\mathcal{F}|$.\medskip

\textit{Case 2: $n \leq 2r-3$.} Then $r > \frac{n}{2} + 1$. 
Let $\mathcal{E}_0' = \Delta_{n,2} \circ \dots \circ 
\Delta_{1,2}(\mathcal{E}_0)$. Thus, $\mathcal{E}_0' \subseteq {\mathcal{L}_{n,2}}^{(r)}$. By Lemma~\ref{complemma}, 
\begin{equation}
\mbox{$E \cap F \cap X_n \neq \emptyset$ for any $E, F \in \mathcal{E}_0'$.} \label{proof.1}
\end{equation} 
For any $E \in \mathcal{E}_0$ and any $F \in \mathcal{E}_1$, we have $\emptyset \neq E \cap F \subseteq X_n$ and $E \cap X_n \subseteq \delta_{i,2}(E)$ for any $i \in [n]$. It clearly follows that
\begin{equation}
\mbox{$E \cap F \cap X_n \neq \emptyset$ for any $E \in \mathcal{E}_0'$ and any $F \in \mathcal{E}_1'$.} \label{proof.2}
\end{equation} 
Let $\mathcal{A} = \{E \cap X_n \colon E \in \mathcal{E}_0'\}$ and $\mathcal{B} = \mathcal{E}_1'$. By (\ref{proof.1}), $\mathcal{A}$ is intersecting. By (\ref{proof.2}), $\mathcal{A}$ and $\mathcal{B}$ are cross-intersecting. 
Let $a_i = {n-i \choose r-i}$ for each $i \in \{0\} \cup [r]$, and let $a_i = 0$ for each $i \in [r+1,n]$. Since $r \leq n-1$, we have
\begin{equation} a_0 > \dots > a_r > a_{r+1} = \dots = a_n = 0. \label{proof.3}
\end{equation}
Let $b_{r-1} = 1$ and let $b_i = 0$ for each $i \in (\{0\} \cup [n]) \backslash \{r-1\}$. 

Consider any $i \leq \frac{n}{2}$. Then $i < r-1$, so $b_i = 0 \leq a_{n-i}$. If $i < n/2$, then $a_i > a_{n-i}$ (by (\ref{proof.3}) as $i < \frac{n}{2} < r$), so $a_i + b_i \geq a_{n-i} + b_{n-i}$. If $i = n/2$, then $n-i = i$, so $a_i + b_i = a_{n-i} + b_{n-i}$. 
(Note that if we ignore Case 1 and rely solely on Theorem~\ref{FJTresult}, then we need to prove Conjecture~\ref{FJTconj} for $n \leq 2r-2$, that is, $r \geq \frac{n}{2}+1$. In this case, we may have $i = r-1$. Suppose $i = r-1$. Since $i \leq \frac{n}{2} \leq r-1 = i$, $i = \frac{n}{2}$. Thus, $b_i = 1 = b_{n-i}$ and $a_i = a_{n-i} \geq 1$, and hence again $a_i + b_i = a_{n-i} + b_{n-i}$ and $a_{n-i} \geq b_i$.)

We have shown that $a_i + b_i \geq a_{n-i} + b_{n-i}$ and $a_{n-i} \geq b_i$ for each $i \leq n/2$. By Theorem~\ref{mainresult}, 
\begin{equation} \sum_{A \in \mathcal{A}} a_{|A|} + \sum_{B \in \mathcal{B}} b_{|B|} \leq \sum_{A \in \mathcal{X}_n} a_{|A|} + \sum_{B \in \mathcal{X}_n} b_{|B|}, \nonumber
\end{equation}
where $\mathcal{X}_n = \{A \subseteq X_n \colon x_1 \in A\}$. Now
\begin{align} |\mathcal{E}_0| &= |\mathcal{E}_0'| = \sum_{i=0}^r \left| \left\{E \in \mathcal{E}_0' \colon E \cap X_n \in \mathcal{A}^{(i)} \right\} \right| \nonumber \\
&\leq \sum_{i=0}^r \left| \left\{E \in {\mathcal{L}_{n,2}}^{(r)} \colon E \cap X_n \in \mathcal{A}^{(i)} \right\} \right| = \sum_{i=0}^r |\mathcal{A}^{(i)}| a_i = \sum_{A \in \mathcal{A}} a_{|A|} \nonumber 
\end{align}
and $|\mathcal{E}_1| = |\mathcal{B}| = \sum_{B \in \mathcal{B}} b_{r-1} = \sum_{B \in \mathcal{B}} b_{|B|}$. Clearly, $|\mathcal{F}_0| = \sum_{A \in \mathcal{X}_n} a_{|A|}$ and $|\mathcal{F}_1| = \sum_{B \in \mathcal{X}_n} b_{|B|}$. We have
\[|\mathcal{E}| = |\mathcal{E}_0| + |\mathcal{E}_1| \leq \sum_{A \in \mathcal{A}} a_{|A|} + \sum_{B \in \mathcal{B}} b_{|B|} \leq \sum_{A \in \mathcal{X}_n} a_{|A|} + \sum_{B \in \mathcal{X}_n} b_{|B|} = |\mathcal{F}_0| + |\mathcal{F}_1| = |\mathcal{F}|.\]
Hence the result.~\hfill{$\Box$}

\end{document}